\documentclass{article}
\usepackage{float}
\usepackage[utf8]{inputenc}
\usepackage{amssymb}
\usepackage{amsmath}
\usepackage{amsfonts}
\usepackage{graphicx}
\usepackage{amsthm}

\begin{document}

\title{The Diamond Theorem}
\author{Steven H. Cullinane}
\date{August 2013}

\maketitle

Finite projective geometry underlies the structure of the 35 square patterns in R. T. Curtis's Miracle Octad Generator, and also explains the surprising symmetry properties of some simple graphic designs-- found, for instance, in quilts.

%%%%%%

\begin{figure}[H]
\begin{center}
\includegraphics{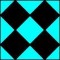}
\end{center}
\caption{Four-diamond figure D made up of 16 tiles in a 4x4 array.}
\end{figure}

We regard the four-diamond figure D above as a 4x4 array of two-color diagonally-divided square tiles.

Let G be the group of 322,560 permutations of these 16 tiles generated by arbitrarily mixing random permutations of rows and of columns with random permutations of the four 2x2 quadrants.

\begin{paragraph}{Theorem} Every G-image of D (as at right in image below) has some ordinary or color-interchange symmetry.
\end{paragraph}

\begin{figure}[H]
\begin{center}
\includegraphics{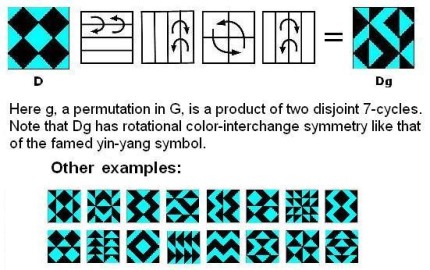}
\end{center}
\caption{Example of graphic permutation.}
\end{figure}

\paragraph{Remarks}

Some of the patterns resulting from the action of G on D have been known for thousands of years. It is perhaps surprising that the patterns' interrelationships and symmetries can be explained fully only by using mathematics discovered just recently (relative to the patterns' age)-- in particular, the theory of automorphism groups of finite geometries.

Using this theory, we can summarize the patterns' properties by saying that G is isomorphic to the affine group A on the linear 4-space over GF(2) and that the 35 structures of the 840 = 35 x 24 G-images of D are isomorphic to the 35 lines in the 3-dimensional projective space over GF(2).

This can be seen by viewing the 35 structures as three-sets of line diagrams, based on the three partitions of the four-set of square two-color tiles into two two-sets, and indicating the locations of these two-sets of tiles within the 4x4 patterns. The lines of the line diagrams may be added in a binary fashion (i.e., 1+1=0). Each three-set of line diagrams sums to zero-- i.e., each diagram in a three-set is the binary sum of the other two diagrams in the set. Thus, the 35 three-sets of line diagrams correspond to the 35 three-point lines of the finite projective 3-space PG(3,2).

For example, here are the line diagrams for the patterns in Fig. 2 above:

\begin{figure}[H]
\begin{center}
\includegraphics{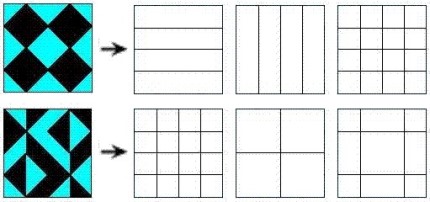}
\end{center}
\caption{Line diagrams indicate patterns' structure.}
\end{figure}

Shown below are the 15 possible line diagrams resulting from row/column/quadrant permutations. These 15 diagrams may, as noted above, be regarded as the 15 points of the projective 3-space PG(3,2).

\begin{figure}[H]
\begin{center}
\includegraphics{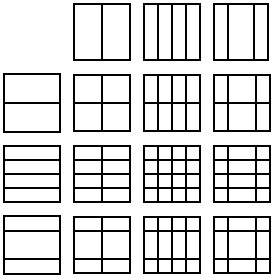}
\end{center}
\caption{Graphic versions of the 15 points of PG(3,2)}
\end{figure}

The symmetry of the line diagrams accounts for the symmetry of the two-color patterns. (A proof shows that a 2nx2n two-color triangular half-squares pattern with such line diagrams must have a 2x2 center with a symmetry, and that this symmetry must be shared by the entire pattern.)

Among the 35 structures of the 840 4x4 arrays of tiles, orthogonality (in the sense of Latin-square orthogonality) corresponds to skewness of lines in the finite projective space PG(3,2).

We can define sums and products so that the G-images of D generate an ideal (1024 patterns characterized by all horizontal or vertical "cuts" being uninterrupted) of a ring of 4096 symmetric patterns. There is an infinite family of such "diamond" rings, isomorphic to rings of matrices over GF(4).

The proof uses a simple, but apparently new, decomposition technique for functions into a finite field.

The underlying geometry of the 4x4 patterns is closely related to the Miracle Octad Generator of R. T. Curtis-- used in the construction of the Steiner system S(5,8,24)-- and hence is also related to the Leech lattice, which, as Walter Feit has remarked, "is a blown up version of S(5,8,24)."

\begin{figure}[H]
\begin{center}
\includegraphics{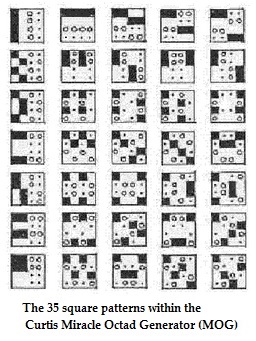}
\end{center}
\caption{The 35 square patterns within the original (1976) version of the Miracle Octad Generator (MOG) of R. T. Curtis}
\end{figure}

As originally presented, the Curtis MOG was a correspondence between the 35 partitions of an 8-set into two 4-sets and the 35 patterns illustrated above. That correspondence was preserved by the actions of the Mathieu group $M_{24}$ on a rectangular array.

The same line diagrams that explain the symmetry of the diamond-theorem figures also explain the symmetry of Curtis's square patterns. The same symmetry group, of order 322,560, underlies both the diamond-theorem figures and the square patterns of the MOG. In the diamond theorem the geometry of the underlying line diagrams shows that this is the group of the affine 4-space over GF(2). In Curtis's 1976 paper this group, under the non-geometric guise $2^{4}.A_{8}$, is shown to be the \emph{octad stabilizer} subgroup of the Mathieu group $M_{24}$.

The above article is an expanded version of Abstract 79T-A37, "Symmetry invariance in a diamond ring," by Steven H. Cullinane, \emph{Notices of the American Mathematical Society}, February 1979, pages A-193, 194.
\vspace{10mm}

\textbf{REFERENCES:}
\vspace{5mm}

R. T. Curtis.  A new combinatorial approach to $M_{24}$. Proceedings of the Cambridge Philosophical Society 79 (1976), 25-42.

\end{document}